\documentclass[12pt]{article}
\usepackage[dvipsnames,usenames]{color}
\usepackage{amsmath}
\usepackage{amsfonts}
\usepackage{amssymb}
\def\d{\delta}
\def\b{\beta}
\def\a{\alpha}
\def\e{\epsilon}
\def\b{\beta}
\def\-{\overline}

\newcommand{\CC}{{\mathbb C}}
\newcommand{\RR}{{\mathbb R}}

\newcommand{\BB}{{\mathbb B}}

\newcommand{\PP}{{\mathbb P}}

\newcommand{\ov}{\overline}
\newcommand{\p}{\partial}
\newcommand{\w}{\widetilde}
\newcommand{\wt}{\widetilde}
\newcommand{\wh}{\widehat}
\newtheorem{theorem}{Theorem}[section]
\newtheorem{lemma}[theorem]{Lemma}

\newtheorem{proposition}[theorem]{Proposition}

\newtheorem{example}[theorem]{Example}

\begin{document}
\author{Xiaojun Huang, Shanyu Ji and Brandon Lee}
\title{\bf CR and Holomorphic Embeddings and Pseudo-conformally Flat Metrics}

\maketitle


\section{Introduction}
In Several Complex Variables , understanding when a CR manifold can
be  embedded into a sphere is a subtle problem. Forstneric
\cite{F86} and Faran \cite{Fa88} proved the existence of real
analytic strictly pseudoconvex hypersurfaces  in $\CC^{n+1}$ which
do not admit any germ of non-constant holomorphic map taking $M$
into sphere $\p\BB^{N+1}$ for any positive integer $N$. Zaitsev
constructed explicit examples for the Forstneric-Faran phenomenon
\cite{Z08}. Meanwhile, there have been much work done to prove the
uniqueness of such embeddings up to the action of automorphisms. For
instance, a well-known rigidity theorem says that that if $M^{2n+1}$
is a CR spherical immersion inside $\p\BB^{N+1}$ with $N\le 2n-1$,
then $M$ must be totally geodesic (i.e., $M$ is the image of
$\p\BB^{n+1}$ by a linear fractional CR map). Ebenfelt, Huang and
Zaitsev (\cite{EHZ04}, Theorem 1.2) proved that if $d<\frac{n}{2}$,
any smooth CR-immersion $f: M \to \p \BB^{n+d+1}$, where $M$ is a
smooth CR hypersurface of dimension $2n+1$, is rigid. Oh in [Oh]
obtained a very interesting result on the non-embeddability for real
hyperboloids into spheres of low codimension.
Kim and Oh \cite{KO06} found a necessary
and sufficient condition for the local holomorphic embeddability
into a sphere of a generic strictly pseudoconvex pseudo-Hermitian CR
manifold in terms of its Chern-Moser curvatures. Along these lines,
we mention recent studies in the papers of Huang-Zhang [HZ],
Ebenfelt-Sun [ES] and Huang-Zaitsev [HZ].
 We also refer the reader to a recent survey paper
[HJ07] by the first two authors and many references therein. Our
fist goal in this paper is to study the non-embddability property
for a class of hypersurfaces, called real hypersurfaces of
involution type, in the low codimensional case, by making use of
property of  a naturally related Gauss curvature. We mention also
the paper by Kolar-Lambel where degenerate revolution hypersurfaces
in ${\CC}^2$ were studied.

\medskip

Consider a {\it real hypersurface of revolution type} defined by

\begin{equation} \label{M}
\begin{split}
& M=\{(z, w)\in {\CC}^{n}\times {\CC}\ | \ r=0\}\\
& r=p(z, \ov z)+q(w, \ov w),\ \ q(w,\-{w})=\ov q(w,\-{w}),\ \
\hbox{d}(q)|_{\{q=0\}}\not =0,\\
&p(z, \ov z)=\sum_{1\le
\a,\b\le n}h_{\alpha\ov\beta} z^\alpha \ov z^\beta.\\
\end{split}
\end{equation}
Here $(h_{\alpha \ov\beta})$ is a positive definite Hermitian
matrix. Such a real hypersurface apparently admits a $U(n)$-action
and was   studied by Webster  in \cite{W02}. Associated with such a
real hypersurface is a domain $D_0$ in ${\CC}$ defined by
$D_0:=\{w\in {\CC}:\ q(w,\-{w})<0\}$. Assume that $M$ is strongly
pseudoconvex in  a certain neighborhood $U_0$ of  $w_0\in D_0$,
Webster observed that then $h:=-(\log q)_{w\-{w}}>0$ in $U_0$ and
thus we have a well-defined Hermitian metric $ds^2=h dw d\-{w}$.
Write the Gauss curvature of such a metric as $K$. Define the Gauss
curvature of this metric by $K=-\frac{1}{h}\frac{\p^2}{\p z\-{\p
z}}\log h$. Write $M_0\subset M$ for an open  piece of $M$ whose
projection to the $w$-space in $U_0$.
We first prove  the following result, which reveals the connection
between the hermitian geometry over $D_0$ and the local smooth CR
embeddability of $M$ into a sphere with lower codimension:

\medskip

\begin{theorem}\label{main thm}
Let $M$ be a strongly pseudoconvex real hypersurface of revolution
in $\CC^{n+1}$ defined as in (\ref{M}) with $2 \le n\le N\le 2n-2$.
Let $D_0$, $U_0$, $K$ and $M_0$ be  just defined as above.  Suppose
the Gauss curvature $K\ge -2$ over $U_0$ and  there is a
non-constant smooth CR map from $M_0$ into $\p\BB^{N+1}$. Then
$K\equiv -2$ over $U_0$ and  the embedding image of $M$ in
$\p\BB^{N+1}$ is totally geodesic, namely, a CR transversal
intersection of an affine complex subspace of dimension $(n+1)$ with
$\p\BB^{N+1}$.
\end{theorem}

\begin{example}
Let $q=|w|^2+\e|w|^4-1$ and $(h_{\a\-\b})=I_{n\times n}$ in
(\ref{M}). Then, for $\e>0$, $M$ admits a non-totally geodesic
holomorphic embedding into the unit sphere in ${\CC}^{n+2}$ through
the map: $(z,w)\mapsto (z,w,\sqrt{\e}w^2)$. However, for $\e<0$, the
Gauss curvature $K$ of $ds^2=-(\log q)_{w\-{w}}dw\otimes d\-{w}$ is
given by $K=-2-4\e+o(1)>-2$ near a neighborhood of $w=0$. (See
Example \ref{curvature}.) Thus, by our theorem and the algebraicity
theorem of the first author [Hu94], $M$ in this setting can not be
locally embedded into $\p\BB^{N+1}$ with $N\le 2n-2$. Hence the
curvature assumption is needed in Theorem \ref{main thm}. Similarly,
let $q=|w|^2+\e|w|^4+|w|^6-1$ with $\e<0, |\e|<<1$. Then $M$ defined
by $r=|z|^2+ |w|^2+\e|w|^4+|w|^6-1=0$  is now compact and strongly
pseudoconvex. Since the Gauss curvature $K$ defined above now is
larger than $-2$ in a neighborhood of $0$ in $D_0$, combing Theorem
\ref{main thm} with the algebraicity theorem of the first author in
[Hu94], we also see that any open piece of $M$ can not be smoothly
CR embedded into $\p\BB^{N+1}$ with $N\le 2n-2$. However, we do not
know if the assumption $N\le 2n-2$ can be dropped.
\end{example}

\medskip
Our  proof of Theorem \ref{main thm} is based on the frame work
established in [EHZ04], computations of Pseudo-Hermitian curvature
tensor  in [We02] and the following rigidity lemma obtained by the
first author:

\bigskip \noindent{\bf Regidity Lemma} [Hu99]:\ \ {\it Let $g_1,
..., , g_k, f_1, ..., f_k$ be holomorphic functions in $z \in\CC^n$
near $0$. Assume $g_j(0) = f_j(0) = 0$ for all $j$. Let $A(z, \ov
z)$ be real-analytic near the origin such that
\begin{equation}
\label{Huang lemma}
\sum^k_{j=1} g_j(z) \ov{f_j(z)} =  |z|^2 A(z, \ov z).
\end{equation}
If $k \le n-1$, then $A(z, \ov z) \equiv 0$ and $\sum^k_{j=1} g_j(z) \ov{f_j(z)}\equiv 0$.}

\bigskip

This rigidity lemma has also played an important role in
understanding many other problems in CR geometry. For instance, the
proof of the third gap theorem \cite{HJY12} is obtained by
repeatedly applying this lemma in subtle ways. In \cite{EHZ04}, a
different formulation of the above lemma was formulated. A new
formulation  of this rigidity lemma is presented in Lemma 2.1 of $\S
2$, and will be used in this paper.

\medskip

Along the same lines of applying the above rigidity lemma, we also
study
rigidity problems for conformal maps between a class of  K\"ahler
manifolds with pseud-conformally flat metrics.
More precisely, we prove the following:

\begin{theorem}
Let $f: (X, \omega) \to (Y, \sigma)$ be a   holomorphic conformal
embedding, where $(X, \omega)$ and $(Y, \sigma)$ are K\"ahler
manifolds with $\dim_\CC X=n$ and $\dim_\CC Y=N$. Suppose $2\le n\le
N\le 2n-1$ and that the curvature tensors of $(X, \omega)$ and $(Y,
\sigma)$ are pseudo-conformally flat. Then  $f(X)$ is a totally
geodesic submanifold of $Y$.
\end{theorem}

\medskip

Here we mention that a holomorphic map $f: (M, \omega) \to (N,
\sigma)$ between Hermitian manifolds $M$ and $N$ is called {\it
conformal} if $f^* \sigma = k \omega$ holds for some  positive
constant $k$ on $M$. A tensor $T_{\alpha \ov\beta \mu\ov\nu}$ over a
complex manifold is called {\it pseudo-conformally flat} (cf.
\cite{EHZ04}) if in any holomorphic chart, we have
\begin{equation}
\label{conformally flat}
T_{\alpha \ov\beta \mu\ov\nu} = H_{\alpha\ov\beta} g_{\mu\ov\nu}
+ \hat H_{\mu\ov\beta} g_{\alpha\ov\nu} + H^*_{\alpha\ov\nu} g_{\mu\ov\beta}
+ \w H_{\mu\ov\nu} g_{\alpha\ov\beta}
\end{equation}
where $(H_{\alpha\ov\beta}), (\hat H_{\alpha\ov\beta}),
(H^*_{\alpha\ov\beta})$  and $(\w H_{\alpha\ov\beta})$ are smoothly
varied Hermitian matrices, and $(g_{\alpha\ov\beta})$ is the
smoothly varied  Hemitian metric, over the chart.

\medskip

Basic examples for Hermitian manifolds with pseudo-conformally flat
curvature tensors are the complex space forms: $\CC^n$ with
Euclidean metric, $\CC\PP^n$ with the Fubini-Study metric and
$\BB^n$ with Poincar\'e metric (see $\S 2$). Other more complicated
examples contain
the Bochner-Kahler manifolds [Br01].

\medskip

Concerning the dimension condition $N\le 2n-1$ in Theorem 1.2, we
recall some related  results on global holomorphic immersions. For
$\CC\PP^n$, Feder proved in 1965 \cite{Fed65} that any holomorphic
immersion $f: \CC\PP^n \to  \CC\PP^N$ with $N\le 2n-1$ has totally
geodesic image (realizing $\CC\PP^n$ as a linear subvariety). For
$X=\BB^n/\Gamma$, Cao and Mok proved in 1990 \cite{CM90} that if $f:
X \to Y$ is a holomorphic immersion where $X$ and $Y$ are complex
hyperbolic space forms of complex dimension $n$ and $N$
respectively, such that $X$ is compact and $N \le  2n- 1$, then $f$
has totally geodesic image. In CR geometry, we have the   rigidity
theorem \cite{Hu99}: if $F: \p\BB^{n+1} \to \p\BB^{N+1}$ is a CR map
which is $C^2$-smooth  with $1\le n \le N\le 2n-1$, then $F$ must be
linear fractional. Also, Mok had constructed an example  [Mok02 ]of
a non-totally geodesic holomorphic isometric embedding from the disc
$\Delta$ into $\Delta^p$. For other related rigidity results, we
refer the reader to the papers by  Calabi \cite{Ca53}, Mok-NG [MN],
Mok [Mok], Yuan-Zhang \cite{YZ12} and many references therein.

\medskip



\section{A tensor version of the rigidity lemma}

We first  reformulate the rigidity lemma mentioned in (\ref{Huang
lemma}) into the following version: (See also related formulations
in [EHZ04])

\begin{lemma}
\label{r lemma} Let $A^{\ a}_{\alpha\beta}$ and $B^{\
a}_{\alpha\beta}$ be complex numbers where $1\le \alpha, \beta\le n,
n+1\le a\le N$. Let $(g_{\alpha \ov\beta})$ and $(G_{a \ov b})$ be
Hermitian matrices with $(g_{\alpha \ov\beta})$ positive definite.
Let $(H^{(l)}_{\alpha \ov\beta}), (\hat H^{(l)}_{\alpha \ov\beta}),
(H^{* (l)}_{\alpha \ov\beta}), (\w H^{(l)}_{\alpha \ov\beta})$ be
Hermitian matrices where $1\le l\le k$. Suppose that $N - n \le n-1$
and that
\begin{equation}\begin{split}
\label{condition} \sum^N_{a, b=n+1} G_{a \ov b} A^{\ a}_{\alpha
\beta} X^\alpha X^\beta \ov{B^{\ b}_{\mu \nu}} \ov{X^\mu} \ov{X^\nu}
= \sum^k_{l=1}(H^{(l)}_{\alpha \ov\beta} g_{\mu \ov\nu} + \hat
H^{(l)}_{\mu \ov\beta} g_{\alpha\ov\nu} +  H^{*(l)}_{\alpha \ov\nu}
g_{\mu \ov\beta} + \w H^{(l)}_{\mu \ov\nu} g_{\alpha
\ov\beta})X^\alpha \-{X^\beta} {X^\mu} \ov{X^\nu}
\end{split}\end{equation}
holds for any $X=(X^\alpha)=(X^\beta)=(X^\mu)=(X^\nu)\in \CC^n$.  Then
\begin{equation}
\label{condition2} \sum^{N}_{a, b=n+1} G_{a \ov b} A^{\
a}_{\alpha\ov\beta} X^\alpha X^\beta \ov{B^{\ b}_{\mu \nu}}
\ov{X^\mu} \ov{X^\nu} \equiv 0,\ \ \ \forall X\in \CC^n.
\end{equation}
\end{lemma}

\medskip

\noindent{\it Proof:}\ \ \ The right-hand-side of (\ref{condition}) is equal to
\begin{equation}
\begin{split}
\label{sum}
& \sum^k_{l=1}(H^{(l)}_{\alpha \ov\beta} g_{\mu \ov\nu} + \hat H^{(l)}_{\mu \ov\beta}
g_{\alpha\ov\nu}
+  H^{*(l)}_{\alpha \ov\nu} g_{\mu \ov\beta} + \w H^{(l)}_{\mu \ov\nu} g_{\alpha
\ov\beta})X^\alpha X^\mu \ov{X^\b} \ov{X^\nu}\\
& = \sum^k_{l=1} \bigg(
 H^{(l)}_{\alpha\ov\beta} X^\alpha \ov{X^\beta} |X|^2
+ \hat H^{(l)}_{\mu \ov\beta} X^\mu\ov{X^\beta}|X|^2
+ H^{*(l)}_{\alpha \ov\nu} X^\alpha \ov{X^\nu} |X|^2
+ \w H^{(l)}_{\mu \ov\nu} X^\mu\ov{X^\nu}|X|^2\bigg) \\
& = |X|^2 \sum^k_{l=1}\bigg( H^{(l)}_{\alpha\ov\beta} X^\alpha \ov{X^\beta} + \hat H^{(l)}_{\mu \ov\beta} X^\mu
\ov{X^\beta} + H^{*(l)}_{\alpha \ov\nu} X^\alpha \ov{X^\nu} + \w H_{\mu \ov\nu} X^\mu \ov{X^\nu}\bigg) =
|X|^2 A(X, \ov X)
\end{split}
\end{equation}
where $A(X, \ov X)$ is some real analytic function of $X$.
Then the left hand side of (\ref{condition}) is equal to
\begin{equation}
\begin{split}
&\sum^{N}_{a,b=n+1} G_{ab}A^a_{\alpha\beta} X^\alpha X^\beta
\ov{B^b_{\mu \nu}} \ov{X^\mu} \ov{X^\nu} = \sum^{N}_{a=n+1} g_a(X)
\ov{h_a(X)}
\end{split}
\end{equation}
where $g_a(X)=\sum_{\alpha, \beta} A^a_{\alpha \beta} X^\alpha
X^\beta$ and $h_a(X)= \sum_{b=n+1}^N\sum_{\alpha, \beta}
\-{G_{ab}}B^b_{\alpha \beta} X^\alpha X^\beta$ are holomorphic
functions. Namely, we have
\[
\sum^{N}_{a=n+1} g_a(X) \ov{h_a(X)} = |X|^2 A(X, \ov X),\ \ \ \forall X\in \CC^n.
\]
By the hypothesis: $N - n < n$, it concludes from (\ref{Huang
lemma}) that $A(X, \ov X)\equiv 0$, and thus (\ref{condition2})
holds. \ \ $\Box$

\medskip

\section{Pseudo-Hermitian geometry}

\noindent{\bf CR submanifold of hypersurface type}\ \ \ \ Let $M$ be
a smooth strictly pseudoconvex $(2n+1)$-dimensional CR submanifold
in $\CC^{n+1}$. We have the complexified tangent bundle $CTM$ which
admits the decomposition $CTM=T^{(1,0)}M\bigoplus T^{(0,1)}M.$ A
non-zero real smooth $1$-form $\theta$ along $M$ is said to be a
contact of $M$ is $\theta|_p$ annihilates $T^{(1,0)}_pM\bigoplus
T^{(0,1)}_pM$ for any $p\in M$. Let $r$ be a local defining function
of $M$. Then $\theta=i \p_z r$ is a contact form of $M$ and any
other contact form is a multiple of $\theta$: $k\theta$ with $k\not
=0$ a smooth function along $M$.

\medskip

Now, fix a contact form $\theta$. Then there is a unique smooth
vector field $T$, called the Reeb vector field such that:  (i)
$\theta (T) \equiv 1$, (ii) $d \theta(T, X) \equiv 0$ for any smooth
tangent vector field $X$ over $M$.
The {\it Levi-form} $L_\theta$ with respect to $\theta$  at $p\in M$
is defined by
\begin{equation}
\label{EHZ2004, 3.2}
L_\theta(u, v):=-id\theta( u \wedge \ov{ v})
=i \theta([u, \ov v]),\ \ \ \forall u,
 v\in T^{1,0}_p(M),\ \forall p\in M.
\end{equation}
Recall that we say $(M, \theta)$ to be  {\it strictly pseudoconvex}
if the Levi-form $L_\theta$ is positive definite for all $z\in M$.

\medskip

Let $T' M$ be the annihilator bundle of ${\cal V}:=T^{(0,1)}M$
which is a rank $n+1$ subbundle of $\CC T^* M$.

\medskip

\noindent{\bf Admissible coframe}\ \ \ \ If we choose a local basis
$L_\alpha$, $\alpha = 1, ..., n$, of $(1,0)$ vector fields (i.e.
sections of $\ov{\cal V}=T^{1,0}_M$), so that $(T, L_\alpha,
L_{\ov\alpha})$ is a frame for $\CC TM:=\CC\otimes TM$ where
$L_{\ov\alpha}=\ov{L_\alpha}$. Then the equation in (ii) above is
equivalent to
\begin{equation}
\label{EHZ2004, 3.3}
d \theta = i g_{\alpha\ov\beta} \theta^\alpha\wedge
\theta^{\ov\beta}.
\end{equation}
Here $\theta^{\ov\beta}=\ov{\theta^\beta}$ and
$(g_{\alpha\ov\beta})$ is the (hemitian) {\it Levi form matrix} and
$(\theta, \theta^\alpha, \theta^{\ov\alpha})$ is the {\it  coframe}
dual to $(T, L_\alpha, L_{\ov\alpha})$. \big(For brevity, we shall
say that $(\theta, \theta^\alpha)$ is the coframe dual to $(T,
L_\alpha)$\big). Note that $\theta$ and $T$ are real whereas
$\theta^\alpha$ and $L_\alpha$ always have non-trivial real and
imaginary parts.

\medskip

Without mentioning $T$, we can complete $\theta$ to a coframe
$(\theta, \theta^\alpha)$ by adding $(1, 0)$-cotangent vectors
(the cotangent vectors that annihilate ${\cal V}$) $\theta^\alpha$. The coframe is called
{\it admissible} if
$\langle\theta^\alpha, T\rangle = 0$, for $\alpha = 1, ..., n$.
As other equivalent definitions, $(\theta, \theta^\alpha)$ is
admissible if (\ref{EHZ2004, 3.3}) holds.


\medskip

\noindent{\bf Pseudo-Hermitian geometry on $M$}\ \ \ \
Observe that (by the uniqueness of the Reeb vector field) for a given contact form
$\theta$ on $M$,
the admissible coframes are determined up to transformations
\[
\w \theta^\alpha = u_\beta^{\ \alpha} \theta^\beta,\ \ (u_\beta^{\ \alpha}) \in GL(\CC^n).
\]
Every choice of a contact form $\theta$ on $M$ is called {\it pseudo-Hermitian structure}
and defines a hemitian metric on ${\cal V}$ (and on $\ov{\cal V}$) via the
(positive-definite)
Levi form (see (\ref{EHZ2004, 3.2})).  For every such $\theta$, Tanaka [T75] and Webster
[W78] defined a
{\it pseudo-Hermitian connection} $\bigtriangledown$ on $\ov{\cal V}$
(and also on $\CC TM$) which is expressed relative to an admissible coframe $(\theta,
\theta^\alpha)$
by
\[
\bigtriangledown L_\alpha = \omega^{\ \beta}_\alpha \otimes L_\beta
\]
where the 1-forms $\omega^{\ \alpha}_\beta$  on $M$ are uniquely determined by the
conditions
\begin{equation}
\label{3.1}
d \theta^\beta = \theta^\alpha \wedge \omega_\alpha^{\ \beta}\ \ mod \wedge (\theta
\wedge \theta^\alpha), \ \
d g_{\alpha \ov\beta} = \omega_{\alpha \ov\beta}+\omega_{\ov \beta \alpha}.
\end{equation}

We may rewrite the first condition in (\ref{3.1}) as
\begin{equation}
\label{3.2}
d \theta^\beta = \theta^\alpha \wedge \omega_\alpha^{\ \beta} + \theta
\wedge \tau^\beta,\ \
\tau^\beta=A^\beta_{\ \ov v} \theta^{\ov v},\ \ A^{\alpha\beta} =
A^{\beta\alpha}
\end{equation}
for a suitably determined torsion matrix $(A^\beta_{\ \ov v})$, where the
last symmetry  relation holds automatically (see [W78]).

\medskip

The {\it pseudo-Hermitian curvature} $R^{\ \beta}_{\alpha\mu\ov{\nu}}$
and $W^{\ \beta}_{\alpha\ \mu}$
of the psuedoHermitian connection is given, in view or [W78, (1.27), (1.41)], by
\begin{equation}
\label{EHZ2004, 3.8}
d \omega_\alpha^{\ \beta} - \omega_\alpha^{\ \gamma} \wedge \omega_\gamma^{\
\beta}
=R^{\ \beta}_{\alpha \ \mu \ov{\nu}} \theta^\mu\wedge \theta^{\ov\nu} + W^{\
\beta}_{\alpha\ \mu} \theta^\mu
\wedge \theta - W^\beta_{\ \alpha \ov\nu} \theta^{\ov\nu}\wedge \theta + i
\theta_\alpha\wedge \tau^\beta - i \tau_\alpha \wedge \theta^\beta.
\end{equation}

\section{Local CR embbedings}

\noindent{\bf Coframes on $f: M \to \hat M$}\ \ \ \
Let $f : M \to \hat M$ be a local CR embedding where $M$ is a strictly pseudoconvex
hypersurface in $\CC^{n+1}$ and
$\hat M$ is a strictly pseudoconvex hypersurface in $\CC^{\hat n+1}$.  We use a $\hat\ $
to denote objects associated to
$\hat M$. We shall also omit the $\hat\ $ over frames and coframes if there is no
ambiguity. It will be clear from the
context if a  form is pulled back to $M$ or not. Under the above assumptions, we identify
$M$ with the submanifold
$f(M)$ and write $M \subset \hat M$. Capital Latin indices $A,B$, etc. will run over
the set $\{1, ... , \hat n\}$. Greek indices $\alpha, \beta$, etc. will run over $\{1,
..., n\}$;
Small Latin indices $a, b$, etc. will run over the complementary set
$\{n + 1, ..., \hat n\}$.

\medskip

Let $(\theta, \theta^\alpha)$ and $(\hat \theta, \hat\theta^A)$ be coframes
on $M$ and $\hat M$ respectively, and recall that $f$ is a {\it CR mapping} if
\[
f^*(\hat \theta) = a\theta,\ \ f^*(\hat\theta^A) =
E^A_{\ \alpha} \theta^\alpha + E^A \theta,
\]
where $a$ is a real-valued function and $E^A_{\ \alpha}$, $E^A$ are complex-valued
functions.
applying $f^*$ to the equation

\medskip

We identify $M$ with the submanifold $f(M)$ of $\hat M$ and write $M
\subset \hat M$. Then the CR bundle ${\cal V}=T^{0,1} M$ is a rank
$n$ subbundle of $\hat {\cal V}=T^{0,1} \hat M$ along $M$. Then
there is a rank $(\hat n-n)$ subbundle $N' M$
consisting of 1-forms on $\hat M$ whose pullbacks to $M$ by $f$
vanish. The subbundle $N' M$ is called the {\it holomorphic conormal
bundle} of $M$ in $\hat M$.

We
write $i^*$ for the standard pull back map and $i_*$ for the
push-forward map. Notice that our consideration is purely local. We
let $p\in M$ and fix a local admissible coframe $\{
\theta,\theta^{\a} \}$ for $M$. Let $T$ be the Reeb vector field
associated with $\theta$. Assume that $\wh{M}$ is a small
neighborhood of $0$ in ${\RR}^{\wh{m}}$, $p=0$ and  $M$ is defined
near $0$ by $x_{j}=0$ with $j=m+1,\cdots,\wh{m}$. First, we can
extend $\theta$ to a contact form of $\wh{M}$ in a neighborhood of
$0$.
Write $x'=(x_1,\cdots,x_m)$.
Define
$\wh{\theta}=u\theta$, with $u(x',0)\equiv 1$. Then
$d\wh{\theta}=du\wedge \theta+ud\theta$. We want $d\wh{\theta}\
\lrcorner T=0$ along $M$.
 For this, we  write $ud\theta\
 \lrcorner T= \sum_{j=1}^{\wh{m}}d_j(x',0) dx_j$. Then, we need to have, along $M$:
 $du=\sum_{j=1}^{\wh{m}}d_j(x',0) dx_j$.  Since $T$ is the Reeb vector
 field for $\theta$ along $M$, we have $d_j(x',0)=0$ for $j\le m$.
 Thus, choose $u=1+\sum_{j=m+1}^{\wh{m}}d_j(x',0)x_j$. Then we have $d\wh{\theta}\
\lrcorner T=0$ along $M$. Now, by the uniqueness of the Reeb vector
field, we see that Reeb vector field $\wh{T}$ of $\wh{\theta}$, when
restricted to $M$, coincides with $T$. Extend $\theta^\a$ to  a
neighborhood of $0$ in $\wh{M}$ to get $\wh{\theta}^\a$, and add
$\wh{\theta}^a$ so that $\{\wh{\theta},\wh{\theta}^\a,
\wh{\theta}^a\}$ forms a basis  for $T'\wh{M}$ near $0$. Apparently,
after a linear change for the forms $\{\wh{\theta}^\a,
\wh{\theta}^a\}$, we can assume that the pull-back of
$\wh{\theta}^a$ to $M$ is zero for each $a=n+1,\cdots, \wh{n}$, the
pull back of $\wh{\theta}^\a$ to $M$ is $\theta^\a$ for
$\a=1,\cdots,n$, $\wh{\theta}$ remains the same, and
$\{\wh{\theta},\wh{\theta}^\a, \wh{\theta}^a\}$ is an admissible
coframe along $\wt{M}$ near $0$.

Next, suppose that $d\theta=\sqrt{-1}g_{\a\-{\b}}\theta^\a\wedge
\theta^{\-\b}$ with $g_{\a\-{\b}}=\d_{\a\b}$ along $M$. We can even
make the Levi form of $\wh{M}$ with respect to the co-frame
$\{\wh{\theta},\wh{\theta}^\a, \wh{\theta}^a\}$ also the identical
matrix along ${M}$. Indeed, let $\{T, L_\a\}$ be the dual frame of
$\{\theta,\theta^\a\}$ along $M$. Extend $L_\a$ to a vector field of
type $(1,0)$ in  a neighborhood of $0$ in $\wh{M}$. Find
$\{\wh{L}_a\}$ so that  $\{ \wt{L}_\a, \wh{L}_a \}$ forms  a base of
vector fields of type $(1,0)$ over $\wh{M}$ with its Levi form along
$\wh{M}$ near $0$ the identical matrix. Let
$\{\wh{\theta},\wh{\theta}^\a, \wh{\theta}^a\}$ be the dual co-frame
of $\{\wh{T},\wh{L}_A\}$. Then along $M$, $<i^*(\wh{\theta}^a),
L_\a>=<\wh{\theta}^a, \wh{L}_\a>|_M=0$; $<i^*(\wh{\theta}^a),
T>=<\wh{\theta}^a, \wh{T}|_M>=0$. Hence the pull back of
$\wh{\theta}^a$ to $M$ is zero. Clearly, the pull-back to
$\wh{\theta}^\a$ to $M$ is ${\theta}^\a$ and
$i^*(\wh{\theta})=\theta$. Assume that
$$d\wh{\theta} =\sqrt{-1} g_{A\-{B}}\wh{\theta}^A\wedge
\wh{\theta}^{\- B}+\sum_{A=1}^{\wh{n}}\big(
e_A(x)\wh{\theta}^A+\-{e_A(x)}\wh{\theta}^{\-A}\big)\wedge
\wh{\theta}.$$
Contracting along $\wh{T}$, we see that $e_A\equiv 0$.
Hence, we see that $\{\wh{\theta},\wh{\theta}^\a, \wh{\theta}^a\}$
is an admissible co-frame.
Now, the Levi form of $\wh{M}$ along $M$  is the identity with
respect to such a frame.

 We  say that the pseudo-Hermitian
structure $(\hat M, \hat \theta)$ is {\it admissible} for the pair
$(M, \hat M)$ if the Reeb vector field $\hat T$ for $\hat \theta$ is
tangent to $M$. With the just obtained co-frame $(\hat \theta,
\hat\theta^A)$ on $\hat M$ where $A=1,2,..., \hat n$, the
holomorphic conormal bundle $N' M$ is spanned by the linear
combinations of the $\hat\theta^a$. Summarizing the above, we see
the following basic fact from [EHZ04]:

\medskip

\begin{proposition} (\cite{EHZ04}, Corollary 4.2)
\label{Prop 4.2} Let $M$ and $\hat M$ be strictly pseudoconvex
CR-manifolds of dimensions $2n+1$ and $2\hat n+1$ respectively. Let
$f : M \to \hat M$ be a CR embedding. If $(\theta, \theta^\alpha)$
is any admissible coframe on $M$, then in a neighborhood of any
point $\hat p \in f(M)$ in $\hat M$ there exists an admissible
coframe $(\hat \theta, \hat \theta^A)$ on $\hat M$ with
$f^*(\hat\theta, \hat \theta^\alpha, \hat\theta^a) =$ $(\theta,
\theta^\alpha, 0)$. In particular, $\hat\theta$ is admissible for
the pair $(f(M), \hat M)$, i.e., the Reeb vector field $\hat T$ is
tangent to $f(M)$. Also, when the Levi form of $M$ with respect to
the co-frame $(\theta, \theta^\alpha)$ is the identical matrix, then
we can also choose $(\hat\theta, \hat\theta^A)$ such that the Levi
form of $\wh{M}$ with respect to  $(\hat \theta, \hat \theta^A)$ is
also the identical matrix.
\end{proposition}

\medskip

If we fix an admissible coframe $(\theta, \theta^\alpha)$ on $M$ and
let $(\hat\theta, \hat\theta^A)$ be an admissible coframe on $\hat
M$ near a point $\hat p \in f(M)$, we shall say $(\hat\theta,
\hat\theta^A)$ is {\it adapted} to $(\theta, \theta^\alpha)$ on $M$
if it satisfies the conclusions of the above Proposition above. We
also normalize the Levi-forms with these frame such that they are
identical.

\medskip

\noindent{\bf Second fundamental form}\ \ \ \
Equation (\ref{3.2}) implies that when $(\theta,
\theta^A)$ is adapted to $M$, if the pseudoconformal connection matrix of $( \hat M ,
\hat \theta)$ is
$\hat \omega^{\ A}_B$ , then that of $(M, \theta)$ is the pullback of $\hat \omega^{\
\alpha}_\beta$. The pulled
back torsion $\hat\tau^\alpha$ is $\tau^\alpha$, so omitting the $\hat\ $
over these pullbacks will not cause any ambiguity and we shall do that from now on. By
the normalization
of the Levi form, the second equation in (\ref{3.1}) reduces to
\begin{equation}
\label{3.3}
\omega_{B \ov A} + \omega_{\ov A B} = 0,
\end{equation}
where as before $\omega_{\ov A B} = \ov{\omega_{A\ov B}}$.

\medskip

The matrix of 1-forms $(\omega^{\ b}_\alpha)$
pulled back to $M$ defines the {\it second fundamental form}
of the embedding $f : M \to \hat M$ . Since $\theta^b = 0$ on $M$, equation
(\ref{3.2}) implies that on $M$,

\begin{equation}
\label{3.4}
\omega^{\ b}_\alpha \wedge \theta^\alpha + \tau^b \wedge \theta = 0 ,
\end{equation}
and this implies that
\begin{equation}
\label{3.5}
\omega^{\ b}_\alpha = \omega^{\ b}_{\alpha \ \beta} \theta^\beta, \
\omega^{\ b}_{\alpha\ \beta}
=\omega^{\ b}_{\beta\ \alpha},\ \ \tau^b=0.
\end{equation}
Following [EHZ04], we identify the CR-normal space $T^{1,0}_p \hat
M/ T^{1, 0}_p M$, also denoted by $N^{1, 0}_p \hat M$ with
$\CC^{\hat n-n}$ by choosing the equivalence classes of $L_a$ as a
basis. Therefore for fixed $\alpha, \beta = 1, ..., n$, we view the
component vector $(\omega^{\ a}_{\alpha\ \beta})_{a=n+1, .., \hat
n}$ as an element of $\CC^{\hat n-n}$. Also view the second
fundamental form as a section over $M$ of the bundle $T^{1, 0}M
\otimes N^{1, 0} \hat M \otimes T^{1, 0}M$,

\section{The Pseudo-conformal geometry}

\noindent{\bf Pseudo-conformal geometry}\ \ \ \
We will need the pseudo-conformal connection and structure equations
introduced by Chern and Moser in [CM74]. Let $Y$ be the bundle of
coframes $(\omega, \omega^\alpha, \omega^{\ov \alpha}, \phi)$ on the real
ray bundle $\pi_E : E \to M$ of all contact
forms defining the same orientation of $M$, such
that $d\omega = i g_{\alpha \ov\beta} \omega^\alpha \wedge
\omega^{\ov\beta}+\omega \wedge \phi$ where $\omega^\alpha \in
\pi^*_E(T'M)$ and $\omega$ is the canonical 1-form on $E$.
In [CM74] it was shown that these forms can be completed to a full set of
invariants on  $Y$ given by the coframe of 1-forms

\begin{equation}
\label{EHZ2004, 3.11}
(\omega, \omega^\alpha, \omega^{\ov\alpha}, \phi, \phi^\alpha_{\ \beta},
\phi^{\ov\alpha}, \psi)
\end{equation}
which define the {\it pseudo-conformal connection} on $Y$.

\begin{equation}\begin{split}
\label{EHZ2004, 3.12}
&\phi_{\alpha\ov\beta}+\phi_{\ov\beta \alpha}=d g_{\alpha\ov\beta},\\
&d\omega=i \omega^\mu\wedge \omega_\mu + \omega \wedge \phi, \\
&d\omega^\alpha=\omega^\mu\wedge \phi_\mu^{\ \alpha} + \omega\wedge
\phi^\alpha,\\
&d\phi=i \omega_{\ov\nu}\wedge \phi^{\ov\nu}+i \phi_{\ov\nu}\wedge
\omega^{\ov\nu}+\omega\wedge \psi,\\
& d \phi^{\ \alpha}_\beta=\phi^{\ \mu}_\beta \wedge \phi^{\ \alpha}_\mu
+ i \omega_\beta\wedge \phi^\alpha - i \phi_\beta\wedge \omega^\alpha - i \delta^{\
\alpha}_\beta
\phi_\mu\wedge \omega^\mu - \frac{\delta^{\ \alpha}_\beta}{2}\psi\wedge \omega + \Phi^{\
\alpha}_\beta,\\
& d\phi^\alpha=\phi\wedge \phi^\alpha + \phi^\mu \wedge \phi^{\ \alpha}_\mu - \frac{1}{2}
\psi
\wedge \omega^\alpha + \Phi^\alpha,\\
& d \psi = \phi\wedge \psi + 2i \phi^\mu \wedge \phi_\mu + \Psi,
\end{split}\end{equation}
where the curvature $2$-forms $\Phi_\beta^{\ \alpha}$, $\Phi^\alpha$ and
$\Psi$ are decomposed as
\begin{equation}\begin{split}
\label{EHZ2004, 3.13}
&\Phi^{\ \alpha}_\beta = S^{\ \alpha}_{\beta\ \mu\ov\nu} \wedge
\omega^{\ov\nu}
+ V_{\beta\ \mu}^{\ \alpha} \omega^\mu \wedge \omega +
V^\alpha_{\ \beta \ov\nu} \omega \wedge \omega^{\ov\nu}, \\
&\Phi^\alpha = V^\alpha_{\ \mu\ov\nu} \omega^\mu \wedge \omega^{\ov\nu}
+P_\mu^{\ \alpha} \omega^\mu\wedge \omega + Q_{\ov\nu}^{\ \alpha}
\omega^{\ov\nu} \wedge \omega, \\
& \Psi=-2i P_{\mu\ov\nu}\omega^{\ov\nu} + R_\mu \omega^\mu\wedge \omega
+R_{\ov\nu} \omega^{\ov\nu}\wedge \omega.
\end{split}\end{equation}
where the functions $S^{\ \alpha}_{\beta\ \mu\ov\nu},
V^{\ \alpha}_{\beta\ \mu}, P^{\ \alpha}_\mu,
Q^{\ \alpha}_{\ov\nu}$ together represent the
{\it pseudo-conformal curvature} of M. \footnote{The indices of
$S^{\ \alpha}_{\beta\ \mu\ov\nu}$ here are interchanged comparing
to [CM74] to make them consistent with indices of
$R^{\ \alpha}_{\beta\ \mu\ov\nu}$ in (\ref{EHZ2004, 3.8}).}
As in [CM74] we restrict our attention here to coframes $(\theta, \theta^\alpha)$ for
which the
Levi form $(g_{\alpha \ov\beta})$ is constant. The $1$-forms $\phi^\alpha,
\phi^{\ov\alpha},
\phi_\beta^{\ \alpha}, \psi$ are uniquely determined
by requiring the coefficients in (\ref{EHZ2004, 3.13}) to satisfy certain symmetry and
trace conditions (see [CM74] and the appendix), e.g.
\[
S_{\alpha\ov\beta\mu\ov\nu} = S_{\mu\ov\beta\alpha\ov\nu}
= S_{\mu\ov\nu\alpha\ov\beta} = S_{\ov\nu\mu\ov\beta\alpha}, \ \
S^{\ \mu}_{\mu\ \alpha\ov\beta} = V^{\ mu}_{\alpha\ \mu} = P^{\ \mu}_\mu = 0.
\]

\medskip

Let us fix a contact form $\theta$ that defines a section $M \to E$. Then any
admissible coframe $(\theta, \theta^\alpha)$ for $T^{1,0}M$ defines a unique section $M
\to Y$ for
which the pullbacks of $(\omega, \omega^\alpha)$ coincide with $(\theta, \theta^\alpha)$
and the pullback
of $\phi$ vanishes. As in [W78], we shall use the same notation for the pulled back forms
on $M$
(that now depend on the choice of the admissible coframe). With this
convention, we have

\begin{equation}
\label{EHZ2004, 3.15}
\theta = \omega, \theta^\alpha = \omega^\alpha, \phi= 0
\end{equation}
on $M$.

 \medskip

\noindent{\bf Relationship between psudo-conformal geometry and pseudo-Hermitian
geometry}\ \ \ \
In view of Webster [W78, (3.8)], the pulled back tangential pseudoconformal curvature
tensor
$S^{\ \beta}_{\alpha\ \mu\ov\nu}$ can be obtained from the tangential pseudo-Hermitian
curvature tensor
$R^{\ \beta}_{\alpha\ \mu\ov\nu}$ in (\ref{EHZ2004, 3.8}) by
\begin{equation}
\label{EHZ2004, (3.16)}
S_{\alpha\ov\beta\mu\ov\nu}=R_{\alpha\ov\beta\mu\ov\nu} - \frac{R_{\alpha\ov\beta}
g_{\mu\ov\nu}
+R_{\mu\ov\beta} g_{\alpha\ov\nu} + R_{\alpha\ov\nu} g_{\mu\ov\beta} +
R_{\mu\ov\nu\alpha\ov\beta}
g_{\alpha\ov\beta}}{n+2}+\frac{R(g_{\alpha\ov\beta} g_{\mu\ov\nu} + g_{\alpha\ov\nu}
g_{\mu\ov\beta})}{(n+1)(n+2)}
\end{equation}
where
\[
R_{\alpha\ov\beta}:=R^{\ \mu}_{\mu\ \alpha\ov\beta}\ \ and\ \ R=R^{\ \mu}_\mu
\]
are respectively the {\it pseudo-Hermitian Ricci and scalar curvature} of $(M, \theta)$.

\medskip

\noindent{\bf Traceless component}\ \ \ \ Following  the termenology
in \cite{EHZ04}, we call a tensor $T_{\alpha_1, ..., \alpha_r,
\ov{\beta_1}, ..., \ov{\beta_s}}^{a_1 ... a_t \ov{b_1} ...
\ov{b_q}}$ {\it pseudo-conformally equivalent to 0} or {\it
pseudo-conformally flat} if it is a linear combination of
$g_{\alpha_i \ov{\beta_j}}$ for $i=1,2,...., r$ and $j=1,2,...,s$.
Two tensors $T_{\alpha \ov\beta \mu\ov\nu}$ and $R_{\alpha \ov\beta
\mu\ov\nu}$ are called {\it conformally equivalent} if  $T_{\alpha
\ov\beta \mu\ov\nu} - R_{\alpha \ov\beta \mu\ov\nu}$ are
pseudo-conformally flat. For any tesnor
$R_{\alpha\ov\beta\mu\ov\nu}$, its {\it traceless component} is the
unique tensor  that is trace zero and that is conformally equivalent
to $R_{\alpha\ov\beta\mu\ov\nu}$. We denote the traceless component
by $[R_{\alpha\ov\beta \mu \ov\nu}]$. Formula (\ref{EHZ2004,
(3.16)}) expresses the fact that $S_{\alpha\ov\beta\mu\ov\nu}$ is
the ``traceless component" of $R_{\alpha\ov\beta\mu\ov\nu}$ (cf.
\cite{EHZ04}, (5.5)):
\begin{equation}
\label{EHZ2004, 5.5}
S_{\alpha\ov\beta\mu\ov{\nu}} = [ R_{\alpha\ov\beta\mu\ov{\nu}} ].
\end{equation}
\medskip

\section{Real Hypersurface of Revolution}

\noindent{\bf Real hypersurfaces of revolution}\ \ \ \ Let
$M=\{(z, w)\ | \ r=0\}$ be a real hypersurface of revolution
in $\CC^{n+1}$ with $n\ge 2$ where
\begin{equation}
\label{def f} r=p(z, \ov z)+q(w, \ov w),\ \ q=\ov q\ \ and\ \ p(z,
\ov z)=h_{\alpha\ov\beta} z^\alpha \ov z^\beta.
\end{equation}
where $(g_{\alpha \ov\beta})$ is a positive definite Hermitian
matrix. Also $d(q)\not =0$ when $q=0$.
\medskip

Define $D:=\{(z, w)\ |\ r<0\}$. As the auxiliary curve and domain in
$\CC$, we define $M_0:=\{w\ |\ q(w, \ov w)=0\}$ and $D_0:=\{(w\ |\
q(w, \ov w)<0\}$.
$M$ is
strictly pseudoconvex if and only if on $D_0:=\{q<0\}$, $h:=-(log\
q)_{w \ov w}=\frac{q_w q_{\ov w} - q q_{w\ov w}}{q^2} >0$. Assume
that $M$ is strictly pseudoconvex. Then $D_0$ admits a Hermitian
metric $d s^2 = h dw d \ov w.$ We denote by $K$ its Gaussian
curvature on $D_0$. It was proved in \cite{W02} that for $w\in D_0$
and $(z, w)\in M$ with $n\ge 2$ and $d q \not=0$, the fourth order
Chern-Moser tensor $S(z, w)=0$ if and only if $K(w)=-2$.

\medskip

\noindent{\bf The pseudo-Hermitian curvature of $M$}\ \ \ \ By
Webster, at the point where $d(q)\not =0$, the pseudo-Hermitian
curvature of $M$ is calculated as
\begin{equation}
\label{curvature} R_{\beta\ov\alpha\rho\ov\sigma} =
-A(g_{\beta\ov\alpha} g_{\rho\ov\sigma} + g_{\rho\ov\alpha} g_{\beta
\ov\sigma}) - B p_\beta p_{\ov\alpha} p_\rho p_{\ov\sigma}
\end{equation}
where
\begin{equation}
\label{A} A=-\frac{Q}{1-Q q},\ g_{\a\-{\b}}=h_{\a\-{\b}}+Q p_\a
p_{\-\b}, \theta=-i\p r, \theta^\a=dz_\a-i\eta^\a\theta,
\eta^\a=g^{\a\-{\b}}\eta_{\-\b},\eta_\a=-Qp_\a;
\end{equation}
and
\begin{equation}
\label{B} B=\frac{Q_{w\ov w}}{q_w q_{\ov w}} + 2 Q
\bigg(\frac{Q_w}{q_w} + \frac{Q_{\ov w}}{q_{\ov w}}\bigg) + 3 Q^3 +
\frac{q|(Q_w / q_w) + Q^2|^2}{1-Q q}
\end{equation}
where $Q=\frac{q_{w \ov w}}{q_w q_{\ov w}}$. Notice that the
formulas above were slightly modified from those in [We02], since we
need $(g_{\a\-\b})$  to be positive definite to apply the
Gauss-Codazzi equation here.

\medskip

Here $B$ can also be calculated as
\begin{equation}
\label{B2}
B = \frac{(K+2)k^2}{q^3 (q_w q_{\ov w})^2}
\end{equation}
where $k=q_w q_{\ov w} - q q_{w\ov w}$. We notice that $B$ is a
real-valued function and $B\le 0$ if and only if $K+2\ge 0$.

\medskip


\noindent{\bf Umbilic points of the fourth order Chern-Moser tensor
$S$}\ \ \ \ Let $S$ be the fourth order Chern-Moser tensor when
$n\ge 2$. (For $n=1$, it is replaced by the Cartan invariant). A
point $(z, w)\in M$ is called a {\it umbilic point} if $S(z, w)=0$.

\medskip

It was proved by Webster \cite{W02} that let $w\in D_0$ and $(z, w)\in M$. Then at
points where $dq\not=0$,
\begin{equation}
\label{Webster}
S(z,w)=0\ \ \text{if and only if}\ \  K(w)=-2.
\end{equation}
If $B\equiv 0$, it implies $K\equiv -2$ by (\ref{B2}).

\medskip

\section{Proof of Theorem 1.1}

Let $M_0$ be a connected open piece of $M=\{(z, w)\ | \ r=0\}$, that
is  strongly pseudoconvex  in $\CC^{n+1}$ with $n\ge 2$. Here $M$ is
as in (\ref{def f}).  Assume  that $M_0$ project down to an open
subset $U_0$ of $D_0$. Suppose that there is a non-constant CR map
$F: M_0 \to \p\BB^{N+1}$. By the Hopf lemma and shrinking $M_0$, we
can assume that $F$ is a CR embedding.  Under the assumption as in
Theorem \ref{main thm},  we then need to prove that $F(M)$ must be
the CR transversal intersection of an affine subspace with the
sphere. After shrinking $M_0$ and thus $U_0$, we can assume that
$q_w\not =0$ over $U_0$.

\medskip

We take an admissible coframe $(\theta, \theta^\alpha)$ on $M$ as
mentioned before with $\theta:=-i \p_z r$ as the contact form.
Fixing any point $p\in M_0$, by Proposition \ref{Prop 4.2}, there
exists a neighborhood $\hat U$ of $\hat p:=F(p)$ in $\p\BB^{N+1}$
and an admissible coframe $(\hat\theta, \hat{\theta^A})$ on $\hat U$
such that $F^*(\hat\theta, \hat{\theta^\alpha}, \hat{\theta^a})
=(\theta, \theta^\alpha, 0)$ on $U$, where $U$ is a neighborhood of
$p$ in $M_0$ such that $F(U)=\hat U$.

\medskip

Consider the pseudo-conformal Gauss equation (cf. (5.9) in \cite{EHZ04})
\begin{equation}
\label{conformal Gauss}
[\hat S(X, X, X, X)] = S(X, X, X, X)+ [\langle II(X, X), \ II(X, X)\rangle],\ \ \
\forall X\in T^{1,0}_{\hat p} F(M),
\end{equation}
where $S$ is the pseudo-conformal curvature of $F(M)$, $\wh S$ is
the restriction of the pseudo-conformal curvature of $\p\BB^{N+1}$
on $F(M)$, and $II(X, X)$ is the second fundamental form of
$F(M)\subset \p\BB^{N+1}$. Here the notation $[\ ]$ in
(\ref{EHZ2004, 5.5}) is used and we can regard $X$ as a vector in
$\CC^n$. Locally it can be written as
\begin{equation}
\label{basic eq}
[\hat S_{\alpha\ov\beta \mu\ov\nu}]= S_{\alpha\ov\beta \mu\ov\nu}
+ [g_{a \ov b}\omega_{\alpha\ \mu}^{\ a} \omega^{\ \ov b}_{\ov\beta\ \ov\nu}]
\end{equation}
where $(\omega^{\ b}_\alpha)$ is the second fundamental form of
$F(M)$ and $\omega^{\ b}_\a =\omega^{\ b}_{\alpha\ \beta}
\theta^\beta$, and $(g_{a \ov b})$ is the (Levi) positive definite
Hermitian matrix. Here $\omega^{\ b}_{\alpha\ \beta}$ are functions
satisfying $\omega^b_{\alpha\ \beta}=\omega^{\ b}_{\beta\ \alpha}$
(cf. \cite{EHZ04}, (4.3) and (5.6)). Recall the facts that the
pseudo-conformal curvature of a sphere vanishes and that we have
\[
S_{\alpha\ov\beta \mu\ov\nu}=[R_{\alpha\ov\beta \mu\ov\nu}]
\]
where $R_{\alpha\ov\beta \mu\ov\nu}$ is the pseudo-Hermitian curvature
induced by the pseudo-Hermitian metric on $F(M)$. Then (\ref{basic eq}) becomes
\begin{equation}
\label{Gauss} 0 = [R_{\alpha\ov\beta \mu\ov\nu}] +  [g_{a \ov b}
\omega_{\alpha\ \mu}^{\ a} \omega^{\ \ov b}_{\ov\beta\ \ov\nu}].
\end{equation}

\medskip

Since $F$ is a local CR embedding, we can identify the
psedudo-Hermitian structure $(M, \theta)$ with $(F(M),
(F^{-1})^*\theta)$. In other words, we can identify the
psedo-Hermitian curvature $R_{\alpha\ov\beta\mu\ov\nu}$ on $F(M)$ as
the pseudo-Hermitian curvature over $M$. Then from
(\ref{curvature}), we have $R_{\alpha\ov\beta\mu\ov\nu} =
-A(g_{\alpha\ov\beta} g_{\mu\ov\nu} + g_{\mu\ov\beta} g_{\alpha
\ov\nu}) - B p_\alpha p_{\ov\beta} p_\rho p_{\ov\nu}$. Since $p(z,
\ov z)=h_{\beta\ov\alpha} z^\beta \ov z^\alpha$, we have
\[
p_\beta=h_{\beta\ov\beta'} \ov z^{\beta'}, \ \ p_{\ov\alpha}=h_{\alpha' \ov\alpha}
z^{\beta'}
\]
and thus
\begin{equation} \label{curvature 2}
\begin{split}
& \sum_{\alpha, \beta, \mu, \nu} p_\alpha p_{\ov\beta} p_\mu p_{\ov\nu}
=\sum_{\alpha, \beta, \mu, \nu, \alpha', \beta', \mu', \nu'} h_{\alpha\ov\alpha'} \ov
z^{\alpha'} h_{\beta'
\ov\beta} z^{\beta'} h_{\mu\ov\mu'} \ov z^{\mu'} h_{\nu' \ov\mu} z^{\nu'}
\\
&=\bigg|\sum_{\beta,\nu, \beta', \nu'}  h_{\beta' \ov\beta}
z^{\beta'} h_{\nu' \ov\nu} z^{\nu'} \bigg|^2
\end{split}\end{equation}

\medskip
Now, as in the proof of lemma \ref{r lemma}, we have the following
computation:
$$A_{\a\-\b}g_{\mu\-\nu}X^\a\-{X^\b}X^\mu\-{X^\nu}=B(X,\-{X})|X|^2,$$
where $|X|^2=g_{\a\-{b}}X^\a\-{X^\b}$ and
$B(X,\-{X})=A_{\a\-{b}}X^\a\-{X^\b}.$ We substitute
(\ref{curvature}) and (\ref{curvature 2}) into (\ref{Gauss}) to
obtain
\begin{equation}\begin{split}
0 & = |X|^2  E(X, \ov X) \\
& - B\ \bigg|\sum_{\beta,\nu, \beta', \nu'}  h_{\beta' \ov\beta}
z^{\beta'} h_{\nu' \ov\nu} z^{\nu'} \ov{X^{\beta}}
\ov{X^{\nu}}\bigg|^2\\
& + \sum_{n+1\le a, b\le N} g_{a \ov b} \omega_{\alpha\ \mu}^{\ a} X^\alpha
X^\mu \omega^{\ \ov b}_{\ov\beta\ \ov\nu} \ov{X^\beta}
\ov{X^\nu},\ \ \ \ \ \ \ \ \ \ \forall X\in \CC^n,\ at\ \hat p\in \hat U
\end{split}\end{equation}
for some real analytic function $E(X, \ov X)$. Since $(N-n)+1 \le
\big((2n-2)-n\big)+1=n-1$, we apply Lemma 2.1 to yield that
\begin{equation}
\label{last id}\begin{split} & - B\ \bigg(\sum_{\beta,\nu, \beta',
\nu'}  h_{\beta' \ov\beta} z^{\beta'} h_{\nu' \ov\nu} z^{\nu'}
\ov{X^{\beta}} \ov{X^{\nu}}\bigg) \bigg(  \sum_{\alpha, \mu,
\alpha', \mu'} \ov{h_{\alpha\ov\alpha'}} \ov{z^{\alpha'}}
\ov{h_{\mu\ov\mu'}}
\ov{z^{\mu'}} X^{\alpha} X^{\mu} \bigg) \\
& + \sum_{a,b=n+1}^{N} g_{a \ov b} \omega_{\alpha\ \mu}^{\ a}
X^\alpha X^\mu \omega^{\ \ov b}_{\ov\beta\ \ov\nu} \ov{X^\beta}
\ov{X^\nu} =0,\ \ \ \ \forall X\in \CC^{n}.
\end{split}
\end{equation}

When $B\le 0$, then both terms in the left hand side of the above
equation is nonnegative. Hence, we get that $B\equiv 0$ over $U_0$
and
\[
\sum^N_{a, b=n+1} g_{a \ov b} \bigg(\omega_{\alpha\ \mu}^{\ a}
X^\alpha X^\mu \bigg) \bigg(\omega^{\ \ov b}_{\ov\beta\ \ov\nu}
\ov{X^\beta} \ov{X^\nu}\bigg)\equiv 0,\ \ \ \forall X\in \CC^n.
\]
Since $(g_{\alpha\ov\beta})$ is Hermitian and positive definite, it
implies $\omega^a_{\alpha\ \mu}=0$, $\forall a, \alpha, \mu$ so that
the second fundamental form of $F(M)$ is zero.

\medskip

Then either by the result of Webster in (\ref{Webster}) or by the
result in \cite{JY10}, $F(M)$ and $M$  must be spherical. Thus
$F(M)$ is in the image $G(\p\BB^{n+1})$ for some linear fractional
map $G: \p\BB^{n+1} \to M\subset \p\BB^{N+1}$, by the well-known
rigidity  result in [Hu99].
The proof of
Theorem 1.1 is complete. \ \ $\Box$

\medskip
\begin{example}\label{curvature}
Let $q=|w|^2+\e|w|^4+\phi(w,\-{w})-1$ with $\e\in {\RR}$ and
$\phi=o(|w|^4)$ being smoothly real-valued. Now $D_0=\{w\in {\CC}:
q<0\}$. $ds^2=-(\log q)_{w\-{w}}dw\otimes d\-{w}$ defines a
Hermitian metric in a neighborhood  of $0\in D_0$. The formula for
its Gauss curvature was derived in [(15), We02]:
$$K=-2+q^3k^{-3}\bigg(kq_{ww\-{w}\-{w}}+q|q_{ww\-{w}}|^2-2\Re(q_{ww\-{w}}q_{\-{w}\-{w}}q_w)+q_{w\-{w}}|q_{ww}|^2\bigg)$$
with $k=q_wq_{\- w}-qq_{w\-{w}}$. By a direct computation, one sees
that $K=-2-4\e+o(|w|)$. Hence, for $\e<0$, we have $K>2$ in a small
neighborhood of $0$ in $D_0$
\end{example}

\medskip


\medskip

\section{Examples of pseudo-conformally flat K\"ahler manifolds}

\noindent{\bf Complex space forms}\ \ \ \
A K\"{a}hler manifold of constant holomorphic sectional curvature is called a {\it
complex space form}.
The  universal complex space forms are $\CC^n, \CC\PP^n$ and $\BB^n$ equippred with the
Kahler metric
\[
h_{ij}=\frac{\delta_{ij}}{1+\kappa|z|^2} - \frac{\kappa z_i \ov z_j}{(1+\kappa |z|^2)^2}
\]
with $\kappa=0, 1$ and $-1$ respectively. Also, $z\in {\CC}^n$ in
the ${\CC}^n$ and ${\PP}^n$ (locally chart in this setting) case;
and $|z|<1$ in the hyperbolic space case.
The curvature tensor is given by
\[
\Theta_{ij}=\kappa\bigg(\sum^n_{k, l=1} h_{k\ov l} d z_k \wedge d
\ov{z_l} \bigg) \delta_{ij} - \kappa \sum^n_{l=1} h_{i \ov l} d
\ov{z_l} \wedge d z_j
\]
and
\[
R_{i\ov j k \ov l}=\kappa(h_{i\ov j} h_{k\ov l} + h_{k\ov j} h_{i\ov
l}),
\]


Complex space forms are certainly pseudo-conformally flat.

\medskip

\noindent{\bf Bochner-K\"ahler manifolds} \ \ \ Let $(M,\omega)$ be
a K\"ahler manifold. Write
$\omega=\sum_{i\-{j}}g_{i\-{j}}dz_i\otimes\-{dz_j}$  in a local
holomorphic chart. The Bochner curvature tensor of  $(M,\omega)$ is
defined as the following tensor:
\[
B_{\beta\ov\alpha\rho\ov\sigma} = R_{\beta\ov\alpha\rho\ov\sigma}  -
\frac{g_{\beta\ov\alpha} R_{\rho\ov\sigma} + g_{\rho\ov\alpha}
R_{\beta\ov\sigma} + g_{\beta\ov\sigma}  R_{\rho\ov\alpha} +
g_{\rho\ov\sigma}  R_{\beta\ov\alpha}}{n+2} +
\frac{R(g_{\beta\ov\alpha} g_{\rho\ov\sigma} + g_{\rho\ov\alpha}
g_{\beta\ov\sigma})}{(n+1)(n+2)}
\]
where $R_{\a\-{b}\gamma\delta}$ is the curvature tensor of
$(M,\omega)$, $R_{\a\-{\b}}$ is the Ricci tensor and $R$ is the
scaler curvature of $(M,\omega)$. $(M,\omega)$ is called a
Bochner-K\"ahler manifold if its Bochner curvature tensor is
identically zero. There have been extensive studies on
Bochner-K\"ahler manifolds in the literature, for which we refer the
reader to the paper of Bryant ([Br01]). Bochner-K\"ahler manifolds
are apparently pseudo-conformally flat in our definition.
\medskip

\section{The proof of the Theorem 1.2}

To prove Theorem 1.2, for any point $u_0\in X$, let $z=(z_1, ...,
z_n)$ be a holomorphic coordinate system of $f(X)$ at $z_0=f(u_0)$,
and $\hat z = (z_1, ..., z_n, z_{n+1}, ..., z_N)$ an extension of
$(z_1, ..., z_n)$ to a coordinate system of $Y$ at $z_0$.
We shall fix the following convention for indices: $1 \le i, j, ...,
\le N$, $1 \le \alpha, \beta,\mu,\nu,\gamma,\d ... \le n, \ n+1 \le
a, b, A,B, ..., \le N$.

\medskip

Let us denote by $\hat g_{ij}$ the Hermitian metric of $(Y, \sigma)$ and $\hat R_{i\ov j
k \ov l}$ the curvature
tensor of this metric on $Y$. Let us denote by $g_{\alpha\beta}$ the restriction metric of
the metric $\hat g_{ij}$ on $f(X)$ and and $R_{\alpha \ov\beta \gamma \ov\sigma}$ the
curvature tensor of this reduced metric $g_{ij}$ on $f(X)$.

\medskip

By the Gauss equation, we have the following equation of tensors:
\begin{equation}
\label{Cartan equation}
\hat R_{\alpha\ov\beta\gamma\ov\delta}|_{f(X)} -
R_{\alpha\ov\beta\gamma\ov\delta} = h^A_{\alpha \gamma} \ov{h^B_{\beta \sigma}}
\hat g_{A\ov B}
\end{equation}
where
$h^A_{\alpha\gamma}=\hat g^{A\ov j} \frac{\p \hat g_{\alpha \ov j}}{\p z^\gamma}$
is the second fundamental form of $f(X)$ in $Y$.

\medskip
Since $(Y, \sigma)$ is pseudo-conformally flat,  the restriction of the curvature also satisfies
\begin{equation}\begin{split}
\label{R}
\hat R_{\alpha\ov\beta\gamma\ov\delta}|_{f(X)} & =
(G_{\alpha\ov\beta} \hat g_{\mu\ov\nu}
+ \hat G_{\mu\ov\beta} \hat g_{\alpha\ov\nu} + G^*_{\alpha\ov\nu} \hat g_{\mu\ov\beta}
+ \w G_{\mu\ov\nu} \hat g_{\alpha\ov\beta})|_{f(X)}
\end{split}\end{equation}
where $G_{\alpha\ov\beta}, \hat G_{\alpha\ov\nu},  G^*_{\alpha\ov\nu}, \w G_{\mu\ov\nu}$
are some
Hermitian matrices on $f(X)$.

\medskip

Since $(X, \omega)$ is pseudo-conformally flat, so is $(f(X),
(f^{-1})^*(\omega))$. Since $f$ is holomorphic conformal,  we have
$(f^{-1})^* \omega = k \sigma|_{f(X)}$ for a positive constant
$k>0$. By the assumption that $(X,\sigma)$ is pseudo-conformally
flat, we conclude that $(f(X), \sigma|_{f(X)})$ is also
pseudo-conformally flat, and hence the curvature tensor
$R_{\alpha\ov\beta\gamma\ov\delta}$ is conformally flat on $f(X)$
and it can be written as
\begin{equation}\begin{split}
\label{R2}
R_{\alpha\ov\beta\gamma\ov\delta} & =
H_{\alpha\ov\beta} g_{\mu\ov\nu}
+ \hat H_{\mu\ov\beta} g_{\alpha\ov\nu} + H^*_{\alpha\ov\nu} g_{\mu\ov\beta}
+ \w H_{\mu\ov\nu} g_{\alpha\ov\beta}
\end{split}\end{equation}
where $H_{\alpha\ov\beta}, \hat H_{\alpha\ov\nu},  H^*_{\alpha\ov\nu}, \w H_{\mu\ov\nu}$
are some Hermitian matrices on $f(X)$.

\medskip

By (\ref{Cartan equation})(\ref{R}) and (\ref{R2}), we have
\begin{equation}\begin{split}
\label{Cartan equation 2}
& (G_{\alpha\ov\beta} g_{\mu\ov\nu}
+ \hat G_{\mu\ov\beta} \hat g_{\alpha\ov\nu} + G^*_{\alpha\ov\nu} \hat g_{\mu\ov\beta}
+ \w G_{\mu\ov\nu} \hat g_{\alpha\ov\beta})(z_0)X^\alpha \-X^\beta {X^\mu} \ov{X^\nu} \\
& - ( H_{\alpha\ov\beta} g_{\mu\ov\nu}
+ \hat H_{\mu\ov\beta} g_{\alpha\ov\nu} + H^*_{\alpha\ov\nu} g_{\mu\ov\beta}
+ \w H_{\mu\ov\nu} g_{\alpha\ov\beta})(z_0) X^\alpha \-X^\beta {X^\mu} \ov{X^\nu} \\
& = (h^A_{\alpha \mu} \ov{h^B_{\beta \nu}} X^\alpha \-{X^\beta}
{X^\mu} \ov{X^\nu}, \hat g_{A\ov B})(z_0)
\end{split}\end{equation}
for any $X=(X^\alpha)=(X^\beta)=(X^\mu)=(X^\nu)\in \CC^n$.

\medskip

By the same calculation as in (\ref{sum}), the left hand side of
(\ref{Cartan equation 2}) is equal to $|X|^2 A(X, \ov X)$. Since
$N-n\le 2n-1 -n = n-1$, we can apply Lemma 3.1 to conclude
\[
\sum^N_{A, B=n+1} h^A_{\alpha\mu}(z_0) X^\alpha\-{ X^\beta}
\ov{h^B_{\beta \nu}(z_0)} {X^\mu} \ov{X^\nu} \hat g_{A\ov B}(z_0)
=0,\ \ \forall X \in \CC^n.
\]
Since the Hermitian metric $(\hat g_{A \ov B}(z_0))$ is positive
definite, $h^A_{\alpha\mu}(z_0)=0$ for all $\alpha, \mu$, and $A$.
Since this holds for any point $z$ in $X$, we have proved that the
second fundamental form of $f(X)$ is identically zero, and hence
$f(X)$ is totally geodesic in $Y$, proving the Theorem 1.2. \ \
$\Box$

\medskip

\bibliographystyle{amsalpha}

\end{document}